\documentclass{article}
\usepackage{amsmath}
\usepackage{amsthm}
\usepackage{amssymb}
\usepackage{hyperref}
\usepackage{amssymb}

\providecommand{\bysame}%
{\makebox[3em]{\hrulefill}\thinspace}

\newtheorem{theorem}{Theorem}[section]

\newtheorem{proposition}{Proposition}

\theoremstyle{definition}
\newtheorem{definition}[theorem]{Definition}
\newtheorem{remark}{Remark}

\newcommand{\al}{\alpha}
\newcommand{\la}{\lambda}
\newcommand{\Z}{\mathbb{Z}}
\newcommand{\R}{\mathbb{R}}

\begin{document}
\title{Random Strict Partitions and Determinantal Point Processes}
\author{Leonid Petrov}
\maketitle

\begin{abstract}
  We present new examples of determinantal point processes with
  infinitely many particles. The particles live on the half-lattice
  $\{1,2,\dots\}$ or on the open half-line $(0,+\infty)$. The  main
  result is the computation of the correlation kernels. They have
  integrable form and are expressed through the Euler gamma function
  (the lattice case) and the classical Whittaker functions (the
  continuous case). Our processes are obtained via a limit transition
  from a model of random strict partitions introduced by Borodin (1997)
  in connection with the problem of harmonic analysis for projective
  characters of the infinite symmetric group.
\end{abstract}

\section{Introduction}
In this paper 
we present new examples
of determinantal point 
processes
and compute 
their correlation kernels.
About determinantal
processes, e.g., see
\cite{Soshnikov2000,Peres2006}
and the recent survey \cite{Borodin2009}.

\subsection{A model of random strict partitions}\label{subsection:M}
We begin with describing 
a family
of
probability measures
on the set of all strict partitions. 
These measures depend
on two real parameters $\al\in(0,+\infty)$
and $\xi\in(0,1)$.

By a \textit{strict partition}
we mean a partition without equal parts, 
that is, a sequence of any length
of the form $\la=(\la_1>\dots>\la_{\ell(\la)})$,
where $\la_i\in\Z_{>0}:=\left\{ 1,2,\dots \right\}$.
Set $|\la|:=\la_1+\dots+\la_{\ell(\la)}$,
this is the \textit{weight} of the partition
(we agree that the empty partition $\la=\emptyset$ has zero weight).

Let $\mathrm{Pl}_n$ (where $n=0,1,2,\dots$)
denote the \textit{Plancherel measure} on the set 
of strict partitions of weight $n$:
\begin{equation}\label{pln}
  \mathrm{Pl}_n(\la):=
  \frac{2^{n-\ell(\la)}\cdot n!}{(\la_1!\dots\la_{\ell(\la)}!)^2}
  \prod_{1\le i<j\le \ell(\la)}\left( \frac{\la_i-\la_j}{\la_i+\la_j} \right)^2,
  \qquad |\la|=n.
\end{equation}
This is a probability measure 
on $\left\{ \la\colon|\la|=n \right\}$
which is an analogue (in the theory
of projective representations of symmetric groups)
of the well-known Plancherel measure on ordinary partitions.
The Plancherel measure on strict partitions 
was studied in, e.g., \cite{Borodin1997,Ivanov1999,Ivanov2006,Petrov2009}.

A certain mixing procedure 
(called \textit{poissonization})
for the Plancherel measures on ordinary
partitions was considered
in \cite{Baik1999,Borodin2000b}.
This procedure leads to determinantal point processes.
In our situation we define the following
\textit{poissonized Plancherel measure} on strict partitions:
\begin{equation}\label{Pl}
  \mathrm{Pl}_\theta:=\sum_{n=0}^{\infty}\frac{(\theta/2)^ne^{-\theta/2}}{n!}\mathrm{Pl}_n,
  \qquad \theta>0,
\end{equation}
that is, we mix the measures $\mathrm{Pl}_n$
on $\left\{ \la\colon|\la|=n \right\}$
using the Poisson distribution 
on the set $\left\{ 0,1,2,\dots \right\}$ 
of indices $n$. 
As a result we obtain a probability measure
on all strict partitions.
In \cite{Matsumoto2005} it was proved that
the poissonized Plancherel measure on strict partitions
gives rise to a Pfaffian point process.
We improve this result and show that 
this point process is determinantal (\S\ref{subsection:Pl}).

In \cite{Borodin1997}
Borodin has introduced a deformation
$\mathsf{M}_n^{(\al)}$
of the Plancherel measure
$\mathrm{Pl}_n$ depending on a parameter $\al>0$
(in \cite{Borodin1997} this parameter is
denoted by $x$):
\begin{equation*}
  \mathsf{M}_n^{(\al)}(\la)=\mathit{const}_{\al,n}\cdot
  \mathrm{Pl}_n(\la)\cdot\prod_{i=1}^{\ell(\la)}\prod_{j=1}^{\la_i}
  \big(j(j-1)+\al\big),\qquad |\la|=n.
\end{equation*}
Here $\mathit{const}_{\al,n}$ is the normalizing constant.
As explained in \cite{Borodin1997}, 
the deformed measure $\mathsf{M}_n^{(\al)}$
preserves certain important properties 
of the Plancherel measure
$\mathrm{Pl}_n$.
For $n=0,1,\dots$,
the measure $\mathrm{Pl}_n$ is the limit
of $\mathsf{M}_n^{(\al)}$ as $\al\to+\infty$.

Similarly
to the mixing of the Plancherel
measures $\mathrm{Pl}_n$ (\ref{Pl}),
we 
consider a mixing of the
deformed measures $\mathsf{M}_n^{(\al)}$.
But now as the mixing distribution we take the 
negative binomial distribution
$\{ (1-\xi)^{\al/2}\frac{(\al/2)_n}{n!}\xi^n \}$
on nonnegative integers $n$
with parameter $\xi\in(0,1)$
(here
$(a)_k:=a(a+1)\dots(a+k-1)$
is the Pochhammer symbol).
As a result we again obtain a probability measure
on the set of all strict partitions.
This measure also gives rise to a determinantal point process.

It is convenient to switch from the parameter $\al>0$ to 
a new parameter $\nu:=\frac12\sqrt{1-4\al}$.
The parameter $\nu$ can be either a real number 
$0\le\nu<\frac12$ (if $0<\al\le\frac14$),
or a pure imaginary number (if $\al>\frac14$).
All our formulas below
are symmetric with respect to the replacement 
of $\nu$ by $(-\nu)$.

We denote the above mixing of the measures $\mathsf{M}_n^{(\al)}$
by $\mathsf{M}_{\nu,\xi}$.
The poissonized Plancherel measure $\mathrm{Pl}_\theta$ is the limit of $\mathsf{M}_{\nu,\xi}$
as $\xi\searrow0$, $\al=\frac14-\nu^2\to+\infty$ such that $\al\xi\to\theta$.
In the sequel we call this limit transition
the \textit{Plancherel degeneration}.

\subsection{Point processes}\label{subsection:pp}
We identify 
every strict partition 
$\la$ with the 
point configuration
$\left\{ \la_1,\dots,\la_{\ell(\la)} \right\}$
on the lattice $\Z_{>0}$.
In this way, 
our two-parameter measure $\mathsf{M}_{\nu,\xi}$
on all strict partitions 
gives rise to a point process $\mathbf{P}_{\nu,\xi}$ on $\Z_{>0}$.
The poissonized Plancherel measure
$\mathrm{Pl}_\theta$
also defines a point process 
on $\Z_{>0}$, denote this process by $\mathbf{P}_{\theta}$.
By the very definition, 
$\mathbf{P}_{\nu,\xi}$ and $\mathbf{P}_{\theta}$
are supported by finite configurations.

The point processes $\mathbf{P}_{\nu,\xi}$ and $\mathbf{P}_{\theta}$
have a general structure described in the
following definition.
\begin{definition}
  Let $\mathbf{P}^{(\psi)}$ be the point
  process on $\Z_{>0}$
  that lives on finite configurations
  and
  assigns the following probability to
  every configuration
  $X=\left\{ x_1,\dots,x_N \right\}$:
  \begin{equation}\label{P}
    \mathbf{P}^{(\psi)}(X):=\mathit{const}\cdot
    \left( U(X) \right)^{2}\cdot\prod_{i=1}^{N}\psi(x_i).
  \end{equation}
  Here $\psi$ is a nonnegative function 
  such that 
  $\sum_{x=1}^{\infty}\psi(x)<\infty$,
  $\mathit{const}$ is the normalizing
  constant and 
  \begin{equation*}
    U(X):=\prod_{1\le i<j\le N}\frac{x_i-x_j}{x_i+x_j}.
  \end{equation*}
\end{definition}

The process $\mathbf{P}_{\nu,\xi}$
has the form (\ref{P})
if as the function $\psi(x)$ we take
\begin{equation}\label{psinuxi}
  \psi_{\nu,\xi}(x):=\frac{\xi^x\cos(\pi\nu)}{2\pi}
  \frac{\Gamma(\frac12-\nu+x)\Gamma(\frac12+\nu+x)}{(x!)^2}.
\end{equation}
The process $\mathbf{P}_{\theta}$ also has the form (\ref{P}) if
as $\psi(x)$ we take $\psi_{\theta}(x):=\frac{\theta^x}{2(x!)^2}$
which is 
the Plancherel degeneration
of
$\psi_{\nu,\xi}(x)$.

\subsection{Correlation kernels}
\subsubsection{The pre-limit kernels}
We observe (\S\ref{subsection:L})
that any point process of the form
$\mathbf{P}^{(\psi)}$ (\ref{P}) is determinantal and
explicitly compute correlation kernels 
in the special cases $\mathbf{P}_{\nu,\xi}$ 
(\S\ref{subsection:hhtk}--\ref{subsection:hhtk-proof})
and $\mathbf{P}_{\theta}$ (\S\ref{subsection:Pl}).
The kernel $\mathsf{K}_{\nu,\xi}$ of the process
$\mathbf{P}_{\nu,\xi}$ has integrable form and is expressed
through the Gauss hypergeometric function.
We call $\mathsf{K}_{\nu,\xi}$ the {\em{}hypergeometric--type
kernel\/}.
In \S\ref{subsection:intint}
we present alternative 
double contour 
integral representations for $\mathsf{K}_{\nu,\xi}$.

For any function $\psi$, 
the correlation kernel $\mathsf{K}$ of $\mathbf{P}^{(\psi)}$
is symmetric. 
However, viewed as an operator in the Hilbert space $\ell^2(\Z_{>0})$,
$\mathsf{K}$ is {\em{}not\/} a projection operator as it happens 
in many other (in particular, random matrix)
models with symmetric correlation kernels.

\subsubsection{Limit transitions}\label{subsection:introlimit}
Recall that the process $\mathbf{P}_{\nu,\xi}$
lives on finite configurations on $\Z_{>0}$.
We consider two limit regimes as $\xi\nearrow1$.

In \S\ref{subsection:gamma}
we examine a limit of $\mathbf{P}_{\nu,\xi}$
on the lattice $\Z_{>0}$.
This limit regime corresponds
to studying the asymptotics of 
smallest parts of 
the random strict partition 
distributed according to the measure $\mathsf{M}_{\nu,\xi}$.

In \S\ref{subsection:whittaker}
we consider a scaling limit of $\mathbf{P}_{\nu,\xi}$.
We embed the lattice $\Z_{>0}$
into 
the half-line
$\R_{>0}$, $x\mapsto(1-\xi)x$, where $x\in\Z_{>0}$,
and then pass to the limit as $\xi\nearrow1$.
This limit regime corresponds
to studying the asymptotics of 
scaled largest parts of 
the random strict partition 
distributed according to the measure $\mathsf{M}_{\nu,\xi}$.

The resulting limit point processes 
live on infinite configurations (on $\Z_{>0}$ and $\R_{>0}$, respectively).
One cannot describe the processes
in terms of probabilities of individual configurations.
We use the description in terms of correlation functions.
We show that both limit processes are determinantal and 
explicitly compute their correlation kernels.
The first kernel $\mathsf{K}^{\mathrm{gamma}}_{\nu}$ is 
expressed in terms of the Euler gamma function,
and the second kernel $\mathcal{K}_{\nu}$ is expressed
in terms 
of the Macdonald functions
(they are certain versions of the Bessel functions).
The kernel $\mathcal{K}_{\nu}$ is called the \textit{Macdonald kernel}.
In \S\ref{subsection:whittaker}
we also give an alternative description of $\mathcal{K}_{\nu}$
in terms of a certain Sturm--Liouville operator.
The Macdonald kernel has already appeared in
the recent paper \cite[\S10.2]{Lisovyy2009}
and also in
\cite[\S5]{Olshanski1998} 
in a different context. 

\subsection{Comparison with other models}

\subsubsection{$z$-measures and log-gas systems}
Our determinantal processes
arise from the measures $\mathsf{M}_n^{(\al)}$
on strict partitions 
introduced in \cite{Borodin1997}
which are related
to the problem of harmonic analysis for projective
characters of the infinite symmetric group.
About projective
representations of symmetric
groups, e.g., see \cite{Schur1911,Hoffman1992,Nazarov1992,Ivanov1999}.

Harmonic analysis 
for ordinary characters 
of the infinite symmetric group
leads to the $z$--measures
on ordinary partitions
\cite{Kerov1993,Kerov2004}.
Determinantal processes 
corresponding to the $z$--measures 
are widely studied, e.g., see
\cite{Borodin1998,Baik1999,Borodin2000a,Borodin2000b,Johansson1999,
Borodin2000riemann,Okounkov2001a,Borodin2005,Borodin2006}.
Recently Strahov studied
another 
example of point processes 
of representation--theoretic origin arising from the
$z$--measures with the deformation (Jack) parameter 2
\cite{Strahov-PP,Strahov-repth}.
In that case the point processes are Pfaffian.
The conventional $z$--measures correspond to the Jack parameter 1.
In Remark \ref{rmk:z} we 
compare some of the operators considered in the present
paper and the corresponding objects for the $z$-measures.

On the other hand,
note a similarity of our model (\ref{P})
to lattice log-gas systems
\cite{Forrester}.
The major difference however is that in our
model the pair interaction is directed by the 
factor $\left( U(X) \right)^{2}$
instead of the conventional 
$\left(V(X)\right)^{\beta}$, where
$V(X):=\prod_{1\le i<j\le N}(x_i-x_j)$.
Lattice log--gas 
systems have representation--theoretic
interpretation for $\beta=2$ (the $z$--measures with Jack parameter 1) 
and for $\beta=1$ or $4$ (the deformed $z$--measures
studied in \cite{Strahov-PP,Strahov-repth}).
Our factor $\left( U(X) \right)^{2}$
comes from the structure of the Plancherel measures
$\mathrm{Pl}_n$
(\ref{pln})
on strict partitions
and is specific to the study of 
projective representations of
symmetric groups.

Moreover, the determinantal and Pfaffian processes 
coming from the $z$--measures on ordinary partitions 
are closely related to orthogonal polynomial
ensembles.
Our model seems to lack this property.

\subsubsection{Shifted Schur measure}\label{subsubsection:SSM}

The $z$--measures with Jack parameter 1
are a special case of
the Schur measure on ordinary partitions 
introduced in \cite{Okounkov1999}.
On strict partitions there exists an analogue
of the Schur measure, namely, the 
shifted Schur measure 
introduced in \cite{Tracy2004shifted}.
In \cite[\S4]{Matsumoto2005}
it was pointed out that the poissonized Plancherel measure $\mathrm{Pl}_{\theta}$
can be interpreted as a special case of the shifted Schur measure.
However, it seems that the measures $\mathsf{M}_{\nu,\xi}$
have no such interpretation.
The correlation functions
of the shifted Schur measure were computed 
in \cite{Matsumoto2005}, they are expressed in terms of 
certain Pfaffians. For the poissonized Plancherel measure
these Pfaffians turn into determinants,
see \S\ref{subsection:Pl} below.

There exists another family of (complex-valued)
probability measures on strict partitions
which under certain specializations becomes $\mathsf{M}_{\nu,\xi}$
or $\mathrm{Pl}_{\theta}$. These measures were introduced
by Rains \cite[\S7]{Rains2000}.
We discuss them below in \S\ref{subsection:Rains}.

\subsection{Acknowledgements}
I am very grateful to Grigori Olshanski
for the setting of the problem, permanent attention
and numerous fruitful discussions, and to
Alexei Borodin for very helpful comments on this
work.

\section{Hypergeometric--type kernel}\label{section:hypergeometric--type kernel}
\subsection{The process $\mathbf{P}^{(\psi)}$ as an L--ensemble}\label{subsection:L}

Let $\mathbf{P}^{(\psi)}$ be the point 
process defined by (\ref{P})
with  arbitrary nonnegative function $\psi(x)$
such that $\sum_{x=1}^{\infty}\psi(x)<\infty$.
Let $\mathsf{L}$ be the following
$\Z_{>0}\times\Z_{>0}$ matrix:
\begin{equation}\label{L}
  \mathsf{L}(x,y):=
  \frac{2\sqrt{xy\psi(x)\psi(y)}}{x+y},\qquad x,y\in\Z_{>0}.
\end{equation}
The condition $\sum_{x=1}^{\infty}\psi(x)<\infty$
ensures that the operator in $\ell^2(\Z_{>0})$ 
corresponding to $\mathsf{L}$ is of trace class.
Therefore, the Fredholm determinant $\det(1+\mathsf{L})$
is well defined. 
\begin{proposition}\label{prop:L--ensemble}
  For every finite subset $X=\left\{ x_1,\dots,x_n \right\}\subset \Z_{>0}$,
  we have 
  $\mathbf{P}^{(\psi)}(X)=\frac{\det \mathsf{L}_X}{\det(1+\mathsf{L})}$,
  where by $\mathsf{L}_X$ we denote the submatrix
  $\left[ \mathsf{L}(x_i,x_j) \right]_{i,j=1}^{n}$.
\end{proposition}
This follows from the 
Cauchy determinant identity
\cite[Ch. I, \S4, Ex. 6]{Macdonald1995}.

Proposition \ref{prop:L--ensemble} 
implies that
the random point process
$\mathbf{P}^{(\psi)}$ is an L--ensemble
corresponding to the matrix $\mathsf{L}$ (e.g., see \cite[\S5]{Borodin2009}).
It follows from general 
properties of determinantal point processes
(for example, see \cite[Prop. 2.1]{Borodin2000a})
that the L--ensemble $\mathbf{P}^{(\psi)}$ 
is determinantal, and
its correlation 
kernel has the form 
$\mathsf{K}={\mathsf{L}}(1+\mathsf{L})^{-1}$.
Since $\mathsf{L}$
is symmetric,
the kernel $\mathsf{K}$ is also symmetric.
However, the operator of the form 
$\mathsf{L}(1+\mathsf{L})^{-1}$ in $\ell^2(\Z_{>0})$ 
cannot be a projection operator.

\subsection{Correlation kernel of the process $\mathbf{P}_{\nu,\xi}$}\label{subsection:hhtk}
Here we 
present explicit expressions for the correlation 
kernel $\mathsf{K}_{\nu,\xi}$ of the point process $\mathbf{P}_{\nu,\xi}$
defined by (\ref{P})--(\ref{psinuxi}).
To shorten the notation, set
\begin{equation*}
  \begin{split}
    \textstyle
    \phi_i(x):=&
    \textstyle
    {}_2F_1\left( -\frac12-\nu+i,-\frac12+\nu+i;x+i;\frac{\xi}{\xi-1} \right),\qquad i=0,1,2,\dots;\\
    \textstyle
    \widetilde\phi(x):=&
    \textstyle
    {}_2F_1\left( \frac32+\nu,-\frac12-\nu;x;\frac{\xi}{\xi-1} \right).
  \end{split}
\end{equation*}
where ${}_2F_1$ is the Gauss hypergeometric function.
Since $x\in\Z_{>0}$ and $i\in\Z_{\ge0}$, 
the third parameter of the above hypergeometric functions
is a positive integer,
therefore, $\phi_i(x)$ and $\widetilde\phi(x)$
are well defined.
Also set
\begin{equation}\label{Xi}
  \textstyle
  \Xi(x,y):=\left\{ 
  \Gamma\left( \frac12-\nu+x \right)
  \Gamma\left( \frac12+\nu+x \right)
  \Gamma\left( \frac12-\nu+y \right)
  \Gamma\left( \frac12+\nu+y \right)
  \right\}^{\frac12},
\end{equation}
where $x,y\in\Z_{>0}$.
Note that due to our assumptions on the parameter 
$\nu$ (\S\ref{subsection:M}), the
above expression in the curved brackets is strictly positive.
\begin{theorem}\label{thm:KSum}
  We have
  \begin{equation}\label{KSum}
    \mathsf{K}_{\nu,\xi}(x,y)
    =
    \frac{2\,\Xi\left( x,y \right)\sqrt{xy}}{x+y}
    \sum_{j=0}^{\infty}
    \frac{\xi^{j+(x+y)/2}(1-\xi)^{-2j}
    \phi_{j+1}(x)\phi_{j+1}(y)}{2^{\delta(j)}(x+j)!(y+j)!
    \Gamma(\frac12-\nu-j)\Gamma(\frac12+\nu-j)},
  \end{equation}
  where $\delta(x):=\delta_{x0}$ is the Kronecker delta.
\end{theorem}
\begin{theorem}\label{thm:KPQ}
  We have
  \begin{equation}\label{KPQ}
    \mathsf{K}_{\nu,\xi}(x,y)=
    \frac{\cos(\pi\nu)}{\pi}\frac{\xi^{\frac{x+y}2}\Xi(x,y)}{\sqrt{x!y!(x-1)!(y-1)!}}\cdot
    \frac{\mathsf{A}(x)\mathsf{B}(y)-\mathsf{B}(x)\mathsf{A}(y)}{x^2-y^2},
  \end{equation}
  where $\mathsf{B}(x)=\phi_1(x)$
  and $\mathsf{A}(x)$ can be written in one of the two
  following forms{\rm{}:\/}

  {\rm{}(1)\/} $\mathsf{A}^{(1)}(x):=x\left( 2\phi_0(x)-\phi_1(x) \right)${\rm{};\/}
  
  {\rm{}(2)\/} $\mathsf{A}^{(2)}(x):=\frac{x}{1+\xi}[2\widetilde\phi(x)-(1-\xi)\phi_1(x)]${\rm{}.\/}
\end{theorem}
If $x=y$, formula (\ref{KPQ}) is also true 
when understood according to the L'Hospital's rule.
This agreement 
is also applicable  
to similar  
formulas below.
\begin{remark}\label{rmk:discrete}
  Note that the kernel $\mathsf{K}_{\nu,\xi}$ given by (\ref{KPQ})
  can be viewed as a discrete
  analogue of an integrable operator
  if as variables we take $x^2$ and $y^2$.
  About integrable operators, e.g., see \cite{Its1990,Deift1999}.
  Discrete integrable operators are discussed in 
  \cite{Borodin2000riemann} and
  \cite[\S6]{Borodin2000a}.
\end{remark}
\begin{remark} 
  There is an identity
  \begin{equation}\label{hgf-identity}
    \phi_0(x)=\frac{\widetilde\phi(x)}{1+\xi}
    -\frac{\xi(1+2\nu-x(1-\xi))}{x(1-\xi^2)}\phi_1(x),
  \end{equation}
  which is a combination
  of 2.8(38), 2.8(39) and 2.9(2) in \cite{Erdelyi1953}.
  Therefore,
  $\mathsf{A}^{(2)}(x)=\mathsf{A}^{(1)}(x)+c\,\mathsf{B}(x)$,
  where $c$ does not depend on $x$. 
  Thus, the kernel (\ref{KPQ})
  with $\mathsf{A}^{(1)}$ is identical to
  the one with $\mathsf{A}^{(2)}$.
  
  Furthermore, all our formulas must be symmetric with respect
  to the replacement of $\nu$ by $(-\nu)$.
  Clearly, $\Xi(x,y)$ and
  all the functions $\phi_i(x)$, $i\in\Z_{\ge0}$, possess this property,
  so the kernel (\ref{KSum}) and the kernel  
  (\ref{KPQ}) with $\mathsf{A}^{(1)}$
  do not change under the substitution $\nu\to(-\nu)$.
  The same holds for the kernel (\ref{KPQ}) with $\mathsf{A}^{(2)}$,
  because  
  from (\ref{hgf-identity}) we have
  $\widetilde\phi(x)|_{\nu\to(-\nu)}=
  \widetilde\phi(x)+\frac{\mathrm{\widetilde c}}{x}\phi_1(x)$, where 
  $\widetilde c$ does not depend on $x$.
\end{remark}

\subsection{Scheme of proof of Theorems \ref{thm:KSum} and \ref{thm:KPQ}}\label{subsection:hhtk-proof}
We begin with the argument similar to \cite{Okounkov2001a},
but instead of the infinite wedge space we take the Fock space
with the orthonormal basis
$v_\la=e_{\la_1}\wedge e_{\la_2}\wedge\dots\wedge e_{\la_{\ell(\la)}}$
indexed by all strict partitions (in particular, $v_{\emptyset}=1$).
A similar space is
used in \cite[\S3]{Matsumoto2005} and \cite[\S5.2]{Vuletic2007}.
By calculations in this Fock space 
we first obtain a
{\em{}Pfaffian formula\/} 
for the correlation functions $\rho_{\nu,\xi}$ 
of the point process $\mathbf{P}_{\nu,\xi}$
(and not a determinantal formula as it was in \cite{Okounkov2001a}).
\begin{proposition}\label{prop:Pfaffian}
  There exists 
  a function $\Phi_{\nu,\xi}\colon\left( \Z\setminus\left\{ 0 \right\} \right)^{2}\to\mathbb{C}$
  such that for every finite subset
  $X=\left\{ x_1,\dots,x_n \right\}\subset\Z_{>0}$
  we have
  \begin{equation*}
    \rho_{\nu,\xi}(X)=(-1)^{\sum_{i=1}^{n}x_i}\cdot\mathop{\mathrm{Pf}}\big({\Phi}(X)\big),
  \end{equation*}
  where $\mathop{\mathrm{Pf}}$ means Pfaffian.
  Here $\Phi(X)$ is the $2n\times 2n$ skew-symmetric matrix
  with rows and columns indexed by
  $x_1,\dots,x_n,-x_n,\dots,-x_1$
  such that the $ij$-th 
  element
  of $\Phi(X)$ above the main diagonal is $\Phi_{\nu,\xi}(i,j)$,
  where
  $i$ and $j$
  take values  
  $x_1,\dots,x_n,-x_n,\dots,-x_1$.
\end{proposition} 
Now we explain how 
one can convert 
the above Pfaffian formula 
for the correlation functions of $\mathbf{P}_{\nu,\xi}$
to a determinantal one.
It turns out that 
$\Phi_{\nu,\xi}$
satisfies the following identities
(here
$x,y\in\Z\setminus\left\{ 0 \right\}$):
\begin{itemize}
  \item If $x\ne y$, then $\Phi_{\nu,\xi}(x,-y)=(-1)^{y}\frac{x+y}{x-y}\Phi_{\nu,\xi}(x,y)$.
  \item If $x\ne-y$, then 
    $\Phi_{\nu,\xi}(y,x)=-\Phi_{\nu,\xi}(x,y)$
    and, moreover,
    $\Phi_{\nu,\xi}(-x,-y)=(-1)^{x+y+1}\Phi_{\nu,\xi}(x,y)$.
  \item If $x\ne0$, then $\Phi_{\nu,\xi}(x,-x)+\Phi_{\nu,\xi}(-x,x)=(-1)^{x}$.
\end{itemize}
Fix a finite subset $X=\left\{ x_1,\dots,x_n \right\}\subset\Z_{>0}$.
Define a $2n\times 2n$ matrix
$C_{kl}:=\delta_{kl}+(-1)^{x_{k\wedge l}}\frac{x_{k\wedge l}-x_n}{x_{k\wedge l}+x_n}I_{\{k+l=2n+1\}}$
($k,l=1,\dots,2n$),
where $k\wedge l$ means the minimum of $k$ and $l$, and $I$ means the indicator.
Clearly, $C$ is invertible.
Using the above identities for $\Phi_{\nu,\xi}$, we obtain
\begin{equation*}
  C\Phi(X)C'=\left[
  \begin{array}{cc}
    0&M\\
    -M'&0
  \end{array}
  \right],
\end{equation*}
where $(..)'$ means the matrix transpose and $M$ 
has format $n\times n$.
It follows from properties of Pfaffians that
$\mathop{\mathrm{Pf}}\big(\Phi(X)\big)=(-1)^{n(n-1)/2}(\det C)^{-1}\det M$.
There exist two diagonal $n\times n$
matrices $D_1$ and $D_2$ such that 
$\det (D_1 D_2)=(-1)^{\sum_{i=1}^{n}x_i}(\det C)^{-1}$ and 
$D_1 M^{\circlearrowleft} D_2=\mathsf{K}_{\nu,\xi}(X)=\left[ \mathsf{K}_{\nu,\xi}(x_i,x_j) \right]_{i,j=1}^n$
for some $\Z_{>0}\times \Z_{>0}$ matrix $\mathsf{K}_{\nu,\xi}$.
Here $M^{\circlearrowleft}$ is the matrix
that is obtained from $M$ 
by rotation by 90 degrees counter--clockwise.
Note that 
$\det\left( M^{\circlearrowleft} \right)=(-1)^{n(n-1)/2}\det M$.

Thus, 
$\rho_{\nu,\xi}(X)=(-1)^{\sum_{i=1}^{n}x_i}\mathop{\mathrm{Pf}}\big(\Phi(X)\big)
=\det\mathsf{K}_{\nu,\xi}(X)$,
which means that $\mathsf{K}_{\nu,\xi}$ is the desired correlation kernel.
The kernel $\mathsf{K}_{\nu,\xi}$ is related to $\Phi_{\nu,\xi}$ as follows:
\begin{equation}\label{KPhi}
  \mathsf{K}_{\nu,\xi}(x,y)=\frac{2(-1)^{y}\sqrt{xy}}{x+y}\Phi_{\nu,\xi}(x,-y),\qquad x,y\in\Z_{>0}.
\end{equation}
We obtain 
explicit expressions
for $\Phi_{\nu,\xi}$
in terms of the Gauss hypergeometric
function.
Their form is similar to 
formulas (3.16) and (3.17) in
\cite{Okounkov2001a}.
These expressions for $\Phi_{\nu,\xi}$ together with
relation 
(\ref{KPhi})
imply 
Theorems \ref{thm:KSum} and \ref{thm:KPQ}, respectively.

\begin{remark}\label{rmk:Lnuxi}
  Let $\mathsf{L}_{\nu,\xi}$
  be the operator defined by (\ref{L})
  with $\psi=\psi_{\nu,\xi}$ given by (\ref{psinuxi}).
  Once formula (\ref{KPQ}) for $\mathsf{K}_{\nu,\xi}$ is obtained,
  one can directly check that
  $\mathsf{K}_{\nu,\xi}=\mathsf{L}_{\nu,\xi}(1+\mathsf{L}_{\nu,\xi})^{-1}$.
  Indeed, this relation is equivalent to 
  $\mathsf{K}_{\nu,\xi}+\mathsf{K}_{\nu,\xi}\mathsf{L}_{\nu,\xi}-\mathsf{L}_{\nu,\xi}=0$,
  and the computation of the matrix product $\mathsf{K}_{\nu,\xi}\mathsf{L}_{\nu,\xi}$ 
  mainly reduces to the
  computation of
  sums of the form 
  $\sum_{k=1}^{\infty}
  \frac{\xi^{k}\Gamma(\frac12+\nu+k)\Gamma(\frac12-\nu+k)}{k!(k-1)!}\frac{f(k)}{k+a}$,
  where $a\ne -1,-2,\dots$ is some constant and $f(k)$ is one of the functions
  $k\phi_0(k)$, $k\phi_1(k)$, or $\phi_1(k)$.
  These sums can be computed using 
  Lemma 3.4 in Appendix in \cite{Borodin2000a}.
\end{remark}

\subsection{Double contour integral representations}\label{subsection:intint}
Here we present
two double contour integral
expressions 
for the hypergeometric--type
kernel $\mathsf{K}_{\nu,\xi}$.
Formulas of this type
are useful in certain limit transitions,
e.g., see \cite{Okounkov2002,Borodin2006,Olshanski2009a}.

To obtain double contour integral formulas 
for the correlation kernel $\mathsf{K}_{\nu,\xi}$, 
we write $\mathsf{K}_{\nu,\xi}$
as the sum (\ref{KSum})
and use the contour integral representation for
the hypergeometric function \cite[Lemma 2.2]{Borodin2006}
combined with the identity 
\cite[2.9(2)]{Erdelyi1953}.
To shorten the notation, set
$\mathsf{g}(x):=\frac{\sqrt{\Gamma(\frac12+\nu+x)\Gamma(\frac12-\nu+x)}}{\Gamma(\frac12+\nu+x)}$.
Note that for our values of $\nu$ (see the end of \S\ref{subsection:M}) 
the expression under the square root is positive for all $x\in\Z$.
\begin{proposition}\label{prop:KOk contour 1}
  For all $x,y\in\Z_{>0}$ we have
  \begin{equation*}
    \begin{split}
      \frac{\mathsf{g}(x)}{\mathsf{g}(y)}\mathsf{K}_{\nu,\xi}(x,y)=
      \frac{2\sqrt{xy}}{x+y}&
      \frac1{(2\pi i)^2}
      \oint\limits_{ \left\{ w_1 \right\}}\oint\limits_{ \left\{ w_2 \right\}}
      {\left( 1-w_1\sqrt\xi \right)^{-\frac12+\nu}
      \left( 1-\frac{\sqrt\xi}{w_1} \right)^{\frac12+\nu}}
      \times\\\times&
      {\left( 1-w_2\sqrt\xi \right)^{-\frac12-\nu}
      \left( 1-\frac{\sqrt\xi}{w_2} \right)^{\frac12-\nu}}
      \frac{w_1^{-x}w_2^{-y}}{w_1w_2-1}dw_1dw_2
      \\-
      \frac{\sqrt{xy}}{x+y}
      \frac{1-\xi}{(2\pi i)^2}
      \oint\limits_{ \left\{ w_1 \right\}}
      \oint\limits_{ \left\{ w_2 \right\}}&
      {\left( 1-w_1\sqrt\xi \right)^{-\frac12+\nu}}
      {\left( 1-\frac{\sqrt\xi}{w_1} \right)^{-\frac12+\nu}}
      \times\\\times&
      {\left( 1-w_2\sqrt\xi \right)^{-\frac12-\nu}
      \left( 1-\frac{\sqrt\xi}{w_2} \right)^{-\frac12-\nu}}
      \frac{dw_1dw_2}{w_1^{x+1}w_2^{y+1}}.
    \end{split}
  \end{equation*}
  The contours $\left\{ w_1 \right\}$ and $\left\{ w_2 \right\}$
  go around $0$ and $\sqrt\xi$ in positive direction
  leaving $1/\sqrt\xi$ outside.
  Moreover, in the first integral
  we have to impose an extra condition:
  the contour
  $\{ w_1^{-1} \}$ lies in the interior of the contour
  $\left\{ w_2 \right\}$.
\end{proposition}
\begin{proposition}
  Let $x,y\in\Z_{>0}$. Then
  \begin{equation*}
    \left.
    \begin{split}
      \frac{\mathsf{g}(-y)}{\mathsf{g}(-x)}\mathsf{K}_{\nu,\xi}&(x,y)\\=
      \frac{\sqrt{xy}}{x+y}&
      \frac{1-\xi}{(2\pi i)^2}
      \oint\limits_{ \left\{ w_1 \right\}}
      \oint\limits_{ \left\{ w_2 \right\}}
      {\left( 1-w_1\sqrt\xi \right)^{-\frac12+\nu+x}}
      {\left( 1-\frac{\sqrt\xi}{w_1} \right)^{-\frac12+\nu-x}}
      \times\qquad\\\times
      \Big( 1-&w_2\sqrt\xi \Big)^{-\frac12-\nu+y}
      \left( 1-\frac{\sqrt\xi}{w_2} \right)^{-\frac12-\nu-y}
      w_1^{-x}w_2^{-y}
      \frac{w_1w_2+1}{w_1w_2-1}\cdot
      \frac{dw_1dw_2}{w_1w_2}.
    \end{split}
    \right.
  \end{equation*}
  Here 
  the contours are as in the first integral
  in Proposition {\rm{}\ref{prop:KOk contour 1}\/}.
\end{proposition}

\subsection{Poissonized Plancherel measure}\label{subsection:Pl}
The poissonized Plancherel measure $\mathrm{Pl}_\theta$
defined by (\ref{Pl}) 
gives rise to the point process $\mathbf{P}_{\theta}$ on $\Z_{>0}$, see \S\ref{subsection:pp}.
Denote the L--operator corresponding
to $\mathbf{P}_{\theta}$
by $\mathsf{L}_{\theta}$
(see \S\ref{subsection:L}).
The operator $\mathsf{L}_{\theta}$
is given by (\ref{L}) with $\psi(x)$ replaced by
$\psi_\theta(x)=\frac{\theta^x}{2(x!)^2}$.
\begin{theorem}\label{thm:Pl}
  The point process $\mathbf{P}_{\theta}$ is
  determinantal with
  the correlation kernel
  \begin{equation*}
    \mathsf{K}_{\theta}(x,y)=\frac{\sqrt{xy}}{x^2-y^2}
    \big(
    2\sqrt\theta J_{x-1}J_y-2\sqrt\theta J_{y-1}J_x-(x-y)J_xJ_y
    \big),\qquad x,y\in\Z_{>0}.
  \end{equation*}
  Here $J_k=J_k(2\sqrt\theta)$ is the Bessel function of the first kind.
\end{theorem}
The correlation kernel $\mathsf{K}_{\theta}$ 
is similar to the {\em{}discrete Bessel kernel\/} 
from \cite{Johansson1999} 
and \cite{Borodin2000b} (but note 
the appearance of additional summands
in $\mathsf{K}_{\theta}$).
The kernel $\mathsf{K}_{\theta}$
is obtained 
from
the hypergeometric--type kernel $\mathsf{K}_{\nu,\xi}$
via the Plancherel degeneration
(see the end of \S\ref{subsection:M}).
Moreover,
one can check that $\mathsf{K}_{\theta}=\mathsf{L}_{\theta}(1+\mathsf{L}_{\theta})^{-1}$
using the identities
for the Bessel functions $J_{k}(2\sqrt{\theta})$
from \S2 in \cite{Borodin2000b}.

The poissonized Plancherel measure is a special
case of the shifted Schur measure
introduced and studied in \cite{Tracy2004shifted,Matsumoto2005}.
In \cite[\S3]{Matsumoto2005}
a Pfaffian formula for the correlation functions of the shifted
Schur measure was obtained.
This Pfaffian formula essentially coincides 
with the Plancherel degeneration of the 
formula from Proposition \ref{prop:Pfaffian}.
Therefore as in \S\ref{subsection:hhtk-proof}
the Pfaffian formula from \cite{Matsumoto2005}
turns into a determinantal formula from 
Theorem \ref{thm:Pl} above.

\subsection{Schur--type measure}\label{subsection:Rains}
The measures
$\mathsf{M}_{\nu,\xi}$ and $\mathrm{Pl}_\theta$
defined in \S\ref{subsection:M}
can be included 
in a wider family 
of
(complex-valued) measures on strict partitions
(and, equivalently, on finite point configurations on $\Z_{>0}$).
The latter measures were introduced in \cite[\S7]{Rains2000}.
They
are similar to the Schur measure
on ordinary partitions introduced in \cite{Okounkov1999}
and are defined as
\begin{equation*}
  \mathfrak{M}(\la):=\frac1Z\pi(s_{(\la\mid\la-1)}),
\end{equation*}
where $\la$ is an arbitrary strict partition,
$s_{(\la\mid\la-1)}$ is the Schur function 
indexed by the Young
diagram written in Frobenius notation as
$(\la_1,\dots,\la_{\ell(\la)}\mid\la_1-1,\dots,\la_{\ell(\la)}-1)$ 
(see \cite[Ch. I, \S1]{Macdonald1995}),
and $Z$ is the normalizing constant.
Here $\pi$ is a {\em{}specialization\/}
of the algebra of symmetric functions $\Lambda$
(that is, a multiplicative homomorphism $\pi\colon\Lambda\to\mathbb{C}$)
such that the series $Z=\sum_{\la}\pi(s_{(\la\mid\la-1)})$ converges.
The difference between $\mathfrak{M}$
and the Schur measure 
is that in $\mathfrak{M}$
we have only one Schur function
instead of two functions for the Schur measure.

The probability measure $\mathsf{M}_{\nu,\xi}$
is obtained from $\mathfrak{M}$ if we take the specialization
defined on the Newton power sums as
\begin{equation}\label{specialization}
  \pi(p_k)=
  {\textstyle(\nu-\frac12)}i^{k}\xi^{k/2},\qquad k=1,2,\dots.
\end{equation}
Here $i=\sqrt{-1}$.
Recall that the Newton power sums 
are algebraically independent generators of $\Lambda$.
Though the specialization (\ref{specialization})
is complex-valued, the values $\mathsf{M}_{\nu,\xi}(\la)$
are real positive for all strict partitions $\la$.
The measure $\mathrm{Pl}_\theta$ is obtained in the same way if 
we take the Plancherel degeneration (see the end 
of \S\ref{subsection:M}) of the specialization 
(\ref{specialization}):
\begin{equation*}
  \pi(p_1)=\sqrt\theta,\qquad
  \pi(p_k)=0,\quad k=2,3,\dots.
\end{equation*}

For a wide class of ``admissible'' specializations
Theorem 7.1 in Rains' paper \cite{Rains2000} gives
a determinantal formula for the measure $\mathfrak{M}$
viewed as a complex-valued measure on point configurations on $\Z_{>0}$. 
Denote by $\mathsf{K}_R(x,y)$ the 
correlation kernel 
\cite[(7.6)]{Rains2000}
under the specialization (\ref{specialization}).
In contrast to our kernel $\mathsf{K}_{\nu,\xi}$ 
(\S\ref{subsection:hhtk}),
$\mathsf{K}_R$ 
is not symmetric.
Numerical computations 
suggest the following
relation between 
$\mathsf{K}_R$ and $\mathsf{K}_{\nu,\xi}$.
Fix any $a\in\Z_{>0}$.
Set 
$F(x):=(-1)^{x}\frac{\sqrt{\mathsf{K}_R(a,x)\mathsf{K}_R(x,a)}}{\mathsf{K}_R(x,a)}$ 
(the expression under the square root is real positive).
Then
$\mathsf{K}_{\nu,\xi}(x,y)=\frac{F(x)}{F(y)}\mathsf{K}_R(x,y)$ for all $x,y\in\Z_{>0}$.
This relation between $\mathsf{K}_R$ and $\mathsf{K}_{\nu,\xi}$
is an instance of a so-called ``gauge transformation'' which does not
change the correlation functions and hence the point process.
However, we do not dispose of a 
rigorous proof of the above relation.

\section{Limit transitions}\label{section:Macdonald kernel}
Recall that the measures $\mathbf{P}_{\nu,\xi}$
defined in Introduction
live on finite configurations on $\Z_{>0}$.
As $\xi\nearrow1$,
the probability $\mathbf{P}_{\nu,\xi}(X)$
of every configuration $X\subset \Z_{>0}$
(given by (\ref{P})--(\ref{psinuxi}))
tends to zero.
However, it is possible
to study limits of $\mathbf{P}_{\nu,\xi}$
as $\xi\nearrow1$ in
spaces larger than the
space of finite configurations
in $\Z_{>0}$.
Here we consider two limit
regimes described in \S\ref{subsection:introlimit}.

\subsection{Limit on the lattice}\label{subsection:gamma}

The space of all (possibly infinite) configurations
on $\Z_{>0}$ 
can be identified with 
$\left\{ 0,1 \right\}^{\Z_{>0}}$.
This is a compact space and 
the point process $\mathbf{P}_{\nu,\xi}$
can be viewed as a probability measure on it.
\begin{theorem}
  As $\xi\nearrow1$, 
  there exists a weak limit of the 
  measures $\mathbf{P}_{\nu,\xi}$ 
  on the space $\left\{ 0,1 \right\}^{\Z_{>0}}$.  
  The limit point process on $\Z_{>0}$
  is supported by infinite configurations and 
  is determinantal with the kernel
  \begin{equation*}
    \begin{split}
      \mathsf{K}^{\mathrm{gamma}}_{\nu}(x,y)=&
      \frac{\sqrt{xy}\cdot\mathrm{ctg}(\pi \nu)}{\pi\,\Xi(x,y)}
      \times\\&\times
      \frac{\Gamma(\frac12+\nu+x)\Gamma(\frac12-\nu+y)-
      \Gamma(\frac12+\nu+y)\Gamma(\frac12-\nu+x)}
      {x^2-y^2}.
    \end{split}
  \end{equation*}
  Here $\Xi(x,y)$ is given by {\rm{}(\ref{Xi})\/}.
\end{theorem}
The proof of this theorem uses certain 
asymptotic relations
for the hypergeometric function, cf.
\cite[\S2]{Borodin2005}.

Similar correlation kernels
expressed in terms of the Euler
gamma function have been studied in \cite{Borodin2005,Olshanski2009a}.
However, it seems that
there is no direct link between
our
point process (corresponding to $\mathsf{K}^{\mathrm{gamma}}_{\nu}$)
and processes from \cite{Borodin2005}.

\subsection{Scaling limit and the Macdonald kernel}\label{subsection:whittaker}

Consider embeddings of $\Z_{>0}$ into $\R_{>0}$
depending on our parameter $\xi\in(0,1)$:
$x\mapsto u:=x(1-\xi)\in\R_{>0}$, $x\in\Z_{>0}$.
\begin{theorem}\label{thm:macdonald kernel}
  Under 
  these embeddings, as $\xi\nearrow1$, the point processes $\mathbf{P}_{\nu,\xi}$
  converge to a determinantal point process $\mathcal{P}_{\nu}$ in $\R_{>0}$.
  The correlation kernel $\mathcal{K}_{\nu}$ of $\mathcal{P}_{\nu}$ can be 
  expressed in terms of the Whittaker
  functions 
  {\rm{}(\/}see \cite[\S6.9]{Erdelyi1953} for definition{\rm{}):\/}
  \begin{equation}\label{whittaker}
    \begin{split}
      \mathcal{K}_{\nu}(u,v)\qquad&\\=\frac{\cos(\pi\nu)}{\pi}&
      \frac{2W_{1,\nu}(u)W_{0,\nu}(v)-2W_{1,\nu}(v)W_{0,\nu}(u)-(u-v)W_{0,\nu}(u)W_{0,\nu}(v)}
      {u^2-v^2},
    \end{split}
  \end{equation}
  and also in terms of the   
  Macdonald functions {\rm{}(\/}see \cite[\S7.2.2]{Erdelyi1953} for definition{\rm{}):\/}
  \begin{equation}\label{macdonald}
    \mathcal{K}_{\nu}(u,v)=\frac{\sqrt{uv}\cos(\pi\nu)}{\pi^2}\frac{u K_{\nu+1}(\frac u2)K_{\nu}(\frac v2)-
    vK_{\nu+1}(\frac v2)K_\nu(\frac u2)}{u^2-v^2}.
  \end{equation}
\end{theorem}
For generalities on
point processes and correlation functions on continuous
spaces, e.g., see the survey \cite{Soshnikov2000}.

Formulas (\ref{whittaker})
and (\ref{macdonald})
are proved
using the asymptotics 
\cite[6.8(1)]{Erdelyi1953}
for the hypergeometric
function:
one should 
write the 
kernel $\mathsf{K}_{\nu,\xi}$
using formula
(\ref{KPQ}) with $\mathsf{A}^{(1)}$ and 
$\mathsf{A}^{(2)}$, respectively.
See also Theorem 5.4 in \cite{Borodin2000a}.

The kernel $\mathcal{K}_{\nu}$ is called the \textit{Macdonald kernel}.
Note that $\mathcal{K}_{\nu}$
is an integrable operator
in the variables $u^2$ and $v^2$ (see also Remark \ref{rmk:discrete}).

The Macdonald kernel 
has already appeared in \cite[\S5]{Olshanski1998}
and \cite[\S10.2]{Lisovyy2009}. 
Observe that the kernel \cite[(5.3)]{Olshanski1998}
takes the form (\ref{macdonald})
if we choose the parameters $z_0=\frac14-\frac\nu2$,
$z_0'=\frac14+\frac\nu2$ and change the 
coordinates as
$\xi=\frac{u^2}{16}$, $\eta=\frac{v^2}{16}$.
Note that for our values of $\nu$ (see the end
of \S\ref{subsection:M}) the parameters 
$z_0$ and $z_{0}'$ 
are of principal or complementary series
(e.g., see 
\cite[\S3.7]{Borodin2007} for definition).
\begin{remark}
  There exists
  a simple connection
  between large $n$ limit
  of the measures $\mathsf{M}_n^{(\al)}$ (see \S\ref{subsection:M})
  and $\xi\nearrow1$ limit of the 
  point processes $\mathbf{P}_{\nu,\xi}$.
  These limits are related via 
  the \textit{lifting} construction
  described in \cite[\S5]{Borodin2000a}.
\end{remark}
\begin{remark}
  One can directly
  check that the kernels (\ref{whittaker}) and (\ref{macdonald})
  are the same.
  To do this, one should express
  $W_{\kappa,\mu}(u)$
  and
  $K_\nu(u)$ 
  through the confluent hypergeometric function
  $\Psi(a,c;u)$ and use
  the identity for $\Psi$
  which 
  follows from 
  \cite[6.6(4)--(6)]{Erdelyi1953}:
  \begin{equation*}
    \begin{split}\textstyle
      \Psi\left( -\frac12-\nu,1-2\nu;u \right)-
      \left( \frac14-\nu^2 \right)\Psi\left( \frac32-\nu,1-2\nu;u \right)-&\\
      \textstyle
      -\Psi\left( -\frac12-\nu,-1-2\nu;u \right)+(1+2\nu)
      \Psi&\textstyle\left( \frac12-\nu,1-2\nu;u \right)=0.
    \end{split}
  \end{equation*}
\end{remark}
Consider an integral
operator in $L^2(\R_{>0})$ with the following 
kernel:
\begin{equation}\label{Lnu}
  \mathcal{L}_{\nu}(u,v):=
  \frac{\cos(\pi\nu)}{\pi}\frac{e^{-\frac{u+v}2}}{u+v},\qquad
  u,v\in\R_{>0}.
\end{equation}
The operators $\mathcal{K}_{\nu}$ and $\mathcal{L}_{\nu}$
satisfy the operator 
relation
$\mathcal{K}_{\nu}=\mathcal{L}_{\nu}(1+\mathcal{L}_{\nu})^{-1}$
which is the same as the
relation between the pre-limit 
operators $\mathsf{K}_{\nu,\xi}$ and $\mathsf{L}_{\nu,\xi}$
on the lattice
(\S\ref{section:hypergeometric--type kernel}).
However, the limit process $\mathcal{P}_{\nu}$ cannot be interpreted
as an L--ensemble 
because it has infinite configurations almost surely.

Now let us present 
another description of the
Macdonald kernel $\mathcal{K}_{\nu}$ (\ref{whittaker})--(\ref{macdonald}). Namely, 
we interpret the operator $\mathcal{K}_{\nu}$ 
as a function (in operator calculus sense)
of a Sturm--Liouville differential operator.
Set 
$f_m(u):=\frac1u W_{0,im}(u)$,
where $m\in[0,+\infty)$ is a parameter.
\begin{proposition}\label{prop:diff}
  {\rm{}(1)\/}
  The operator $\mathcal{K}_{\nu}$ commutes
  with the second order differential 
  operator 
  \begin{equation*} 
    \mathfrak{D}=-\frac{d}{du}\,u^2\frac{d}{du}+\frac{1}4u^2.
  \end{equation*}
  That is, the Macdonald kernel $\mathcal{K}_{\nu}(u,v)$
  satisfies
  $\mathfrak{D}_u\mathcal{K}_{\nu}(u,v)=
  \mathfrak{D}_v\mathcal{K}_{\nu}(u,v)$,
  where the subscript $u$ or $v$ indicates the variable
  on which the differential operator acts.

  {\rm{}(2)\/} For every $m\ge0$ we have
  \begin{equation*}
    \mathfrak{D}f_m=\left( m^2+\frac14 \right)f_m,\qquad
    \mathcal{K}_{\nu} f_m=\frac{\cos(\pi\nu)}{\cos(\pi\nu)+\cosh(\pi m)}f_m.
  \end{equation*}
\end{proposition}
One can say that 
$\mathcal{K}_{\nu}=h_\nu(\mathfrak{D})$, where
\begin{equation*}
  h_\nu(r):=\frac{\cos(\pi \nu)}{\cos(\pi\nu)+\cosh\left( \pi\sqrt{r-\frac14} \right)},\qquad
  r\ge\frac14.
\end{equation*}
Note that the operators $\mathfrak{D}$ and $\mathcal{K}_{\nu}$
are self-adjoint in $L^2(\R_{>0})$.

Proposition \ref{prop:diff}
is suggested by results in \cite{Olshanski1998}.
Indeed, observe that the above operator $\mathcal{L}_{\nu}$
coincides with the operator $A$ \cite[(2.29)]{Olshanski1998}
if we set $a=0$ and $\sigma=\cos(\pi\nu)$.
Thus, from \cite[\S3]{Olshanski1998}
it follows that $\mathcal{L}_{\nu}$ commutes with $\mathfrak{D}$ and that 
$\mathcal{L}_{\nu} f_m=\frac{\cos(\pi\nu)}{\cosh(\pi m)}f_m$.
This implies Proposition \ref{prop:diff}
because $\mathcal{K}_{\nu}=\mathcal{L}_{\nu}(1+\mathcal{L}_{\nu})^{-1}$.

The functions $\left\{ f_m \right\}_{m\ge0}$
form a continual basis in $L^2(\R_{>0})$, 
and an explicit Plancherel formula
\cite[(3.5)--(3.6)]{Olshanski1998} (where one must set $a=0$)
holds.
The operators 
$\mathcal{K}_{\nu}$ for our values of $\nu$ 
(see \S\ref{subsection:M})
form a commutative family.
As $\nu\to i\infty$, 
the spectrum of $\mathcal{K}_{\nu}$ becomes closer to $1$,
and the norm of $\mathcal{L}_{\nu}$ tends to infinity.

\begin{remark}\label{rmk:z}
There are certain formal relations between
some of the operators considered in the present
paper and the corresponding objects for the $z$-measures.
More precisely, the 
pairs of corresponding objects are:
the $\Z_{>0}\times\Z_{>0}$ matrix
$\mathsf{L}_{\nu,\xi}$ from Remark \ref{rmk:Lnuxi}
and the matrix \cite[(3.3)]{Borodin2000a};
the operator $\mathcal{L}_\nu$ in $L^2(\R_{>0})$ given by (\ref{Lnu})
and the 
operator \cite[(2.28)--(2.30)]{Olshanski1998};
the Macdonald kernel $\mathcal{K}_\nu$ from Theorem \ref{thm:macdonald kernel}
and the matrix Whittaker kernel from \cite[\S5]{Borodin2000a}.
The relations 
between these objects
are realized by taking \textit{non-admissible}
values of the parameters $(z,z')$ of the $z$-measures, 
namely, $z=-z'=\nu-\frac12$.
Clearly, for our values of $\nu$ the parameters
$(z,z')$ are not of principal or complementary series.
Therefore, it seems that there is no 
direct
connection between
our model and the $z$-measures at the 
level of random point processes.

Let us describe how the operator 
\cite[(2.28)--(2.30)]{Olshanski1998}
is related to $\mathcal{L}_\nu$.
If we set $z=-z'=\nu-\frac12$, the parameter $a$ in 
\cite[(2.29)--(2.30)]{Olshanski1998} vanishes, and 
the parameter $\sigma=\sqrt{\sin(\pi z)\sin(\pi z')}$ should be understood
as $\sigma=i\cos(\pi\nu)$.
We see that the operator 
\cite[(2.28)]{Olshanski1998} 
takes the form
\begin{equation*}
  \left[
  \begin{array}{cc}
    0&i\mathcal{L}_\nu\\
    -i\mathcal{L}_\nu&0
  \end{array}
  \right].
\end{equation*}
This fact implies that 
under the above choice of non-admissible values of 
$(z,z')$ we have\footnote{This formula was suggested
to the author by A.~Borodin in a private communication.}
\begin{equation}\label{KK}
  \mathcal{K}_\nu=\mathcal{K}_{++}-i\mathcal{K}_{+-},
\end{equation}
where $\mathcal{K}_{++}$ and $\mathcal{K}_{+-}$ are the 
blocks of the matrix Whittaker kernel, see 
\cite[\S5]{Borodin2000a}.

The pre-limit relation between 
$\mathsf{L}_{\nu,\xi}$
and the matrix \cite[(3.3)]{Borodin2000a}
has a more complicated structure and 
involves the same non-admissible $(z,z')$.
However, it seems that 
(\ref{KK}) does not have a pre-limit analogue,
that is, there is no tractable relation
between the pre-limit correlation kernel $\mathsf{K}_{\nu,\xi}$
and the matrix hypergeometric kernel from \cite[\S3]{Borodin2000a}.
\end{remark}

\medskip

Dobrushin Mathematics Laboratory,
Kharkevich Institute for Information Transmission Problems,
Bolshoy Karetny per. 19, Moscow, 127994, Russia.

\smallskip

E--mail: \texttt{lenia.petrov@gmail.com}

\end{document}